\documentclass[11pt]{amsart}
\usepackage{amsmath}
\usepackage{amssymb}
\usepackage{latexsym}
\usepackage{graphicx}

\newtheorem{theorem}{Theorem}
\newtheorem{corollary}{Corollary}
\theoremstyle{remark}

\begin{document}

\title[inclusion theorem]{\large A sharp inclusion theorem for a certain class of analytic functions defined by the Salagean derivative}

\author[K. O. Babalola and T. O. Opoola]{K. O. BABALOLA and T. O. Opoola}

\begin{abstract}
In this short note we employ the Briot-Bouquet differential subordination to determine the best possible inclusion relation within a certain family of analytic functions defined by the Salagean derivative.
\end{abstract}



\maketitle

\section{Introduction}
Let $A$ denote the class of functions:
\[
f(z)=z+a_2z^2+...
\]
which are analytic in the unit disk $E=\{z\in \mathbb{C}\colon
|z|<1\}$. Also let P be the class of functions:
\begin{equation}
p(z)=1+p_1z+...\, \label{1}
\end{equation}
which are also analytic in the unit disk $E$ and have positive real part. In \cite{TOO}, Opoola introduced the subclass $T_n^\alpha (\beta)$ consisting of functions $f \in A$ which
satisfy:
\begin{equation}
Re \frac{D^nf(z)^\alpha}{\alpha^nz^\alpha}>\beta\, \label{2}
\end{equation}
\\where $\alpha>0$ is real, $0\le\beta<1$, $D^n(n\in N_0=\{0, 1, 2,
...\})$ is the Salagean derivative operator defined as:
$D^nf(z)=D(D^{n-1}f(z))=z[D^{n-1}f(z)]^{\prime}$ with
$D^0f(z)=f(z)$ and powers in ~(\ref{2}) meaning principal
determinations only. The geometric condition ~(\ref{2}) slightly
modifies the one given originally in \cite{TOO} (see \cite{BO}).

In earlier works \cite{BO, TOO}, inclusion relations have been discussed for the family $T_n^\alpha (\beta)$. In particular, it has been shown that members of the family are related by the inclusion.

\begin{theorem}[\cite{BO, TOO}]
\[
T_{n+1}^\alpha (\beta)\subset T_n^\alpha (\beta), \;\;\;n\in N_0.
\]
\end{theorem}

The object of the present paper is to sharpen the above result. Our result is the following:
\begin{theorem}
\[
T_{n+1}^\alpha (\beta)\subset T_n^\alpha (\delta(\alpha, \beta)), \;\;\;n\in N_0
\]
where
\begin{equation}
\delta(\alpha, \beta)=1+2(1-\beta)\sum_{k=1}^\infty\frac{\alpha}{\alpha+k}(-1)^k.\, \label{3}
\end{equation}
The result is sharp.
\end{theorem}

We will make use of the powerful technique of Briot-Bouquet differential subordination to prove the above result. This technique has been employed frequently in recent times to sharpen and improve many results in geometric function theory. A function $p(z)$ given by ~(\ref{1}) is said to satisfy the Briot-Bouquet differential subordination if 
\begin{equation}
p(z)+\frac{zp^\prime(z)}{\eta p(z)+\gamma}\prec h(z),\;\;\;z\in E\, \label{4}
\end{equation}
where $\eta$ and $\gamma$ are complex constants and $h(z)$ a complex function satisfying $h(0)=1$, and $Re(\eta h(z)+\gamma)>0$ in $E$. It is well known that if $p(z)$ given by ~(\ref{1}) satisfies the Briot-Bouquet differential subordination, then $p(z)\prec h(z)$ \cite{EM}.

A univalent function $q(z)$ is said to be a dominant of ~(\ref{4}) if $p(z)\prec q(z)$ for all $p(z)$ satisfying ~(\ref{4}). If $\widetilde{q}(z)$ is a dominant of ~(\ref{4}) and $\widetilde{q}(z)\prec q(z)$ for all dominants $q(z)$ of ~(\ref{4}), then $\widetilde{q}(z)$ is said to be the best dominant of ~(\ref{4}). The best dominant is unique up to rotation. Furthermore, it is known \cite{EM} that if the differential equation
\begin{equation}
q(z)+\frac{zq^\prime(z)}{\eta q(z)+\gamma}=h(z),\;\;\;q(0)=1
\end{equation}
has univalent solution $q(z)$ in $E$, then $p(z)\prec q(z)\prec h(z)$ and $q(z)$ is the best dominant.

Now let $f\in A$. If we let
\[
p(z)=\frac{D^nf(z)^\alpha}{\alpha^nz^\alpha}
\]
Then we find that
\[
p(z)+\frac{zp^\prime(z)}{\alpha}=\frac{D^{n+1}f(z)^\alpha}{\alpha^{n+1}z^\alpha}
\]
Thus it follows that

{\rm(a)} if $\frac{D^{n+1}f(z)^\alpha}{\alpha^{n+1}z^\alpha}\prec h(z)$, then $\frac{D^nf(z)^\alpha}{\alpha^nz^\alpha}\prec h(z)$.

{\rm(b)} if the differential equation $q(z)+\frac{zq^\prime(z)}{\eta q(z)+\gamma}=h(z)$, $q(0)=1$, has univalent solution $q(z)$ in $E$, then $\frac{D^{n+1}f(z)^\alpha}{\alpha^{n+1}z^\alpha}\prec h(z)$ implies $\frac{D^nf(z)^\alpha}{\alpha^nz^\alpha}\prec q(z)\prec h(z)$ and $q(z)$ is the best dominant.

In view of the above expositions we proceed to the proof of our result in the next section .

\section{Proof of Theorem}
Let $f\in T_{n+1}^\alpha (\beta)$. Define
\[
h(z)=h_\beta(z)=\frac{1+(1-2\beta)z}{1-z}.
\]
It is obvious that $h_\beta(z)$ maps the unit disk onto the plane $Re\;\omega>\beta$. Hence
\[
f\in T_{n+1}^\alpha (\beta)\Leftrightarrow\frac{D^{n+1}f(z)^\alpha}{\alpha^{n+1}z^\alpha}\prec h_\beta(z)
\]
and the differential equation
\[
q(z)+\frac{zq^\prime(z)}{\alpha}=h_\beta(z),\;\;\; q(0)=1
\]
has univalent solution
\[
q_\beta(z)=\frac{\alpha}{z^\alpha}\int_0^zt^{\alpha-1}h_\beta(t)dt.
\]
This yields
\[
q_\beta(z)=1+2(1-\beta)\sum_{k=1}^\infty\frac{\alpha}{\alpha+k}z^k.
\]
That is,
\[
q_\beta(z)=\beta+(1-\beta)\left\{1+2\sum_{k=1}^\infty\frac{\alpha}{\alpha+k}z^k\right\}.
\]

In \cite{BO} we have shown that
\[
Re\left\{1+2\sum_{k=1}^\infty\frac{\alpha}{\alpha+k}z^k\right\}\geq 1+2\sum_{k=1}^\infty\frac{\alpha}{\alpha+k}(-r)^k,\;\;\;|z|=r.
\]
Hence taking limit as $r\rightarrow 1^{-}$, we have  
where $\min_{|z|\leq 1}Re\;q_\beta(z)=q_\beta(-1)=\delta(\alpha,\beta)$ where $\delta(\alpha,\beta)$ is given by ~(\ref{3}). Hence we have
\[
f\in T_{n+1}^\alpha (\beta)\Rightarrow\frac{D^nf(z)^\alpha}{\alpha^nz^\alpha}\prec q_\beta(z)\prec h_\beta(z)
\]
That is,
\[
f\in T_{n+1}^\alpha (\beta)\Rightarrow Re\frac{D^nf(z)^\alpha}{\alpha^nz^\alpha}\geq \delta(\alpha, \beta)
\]
i.e. $f\in T_n^\alpha(\delta(\alpha, \beta))$. The inclusion is sharp since $q_\beta(z)$ is the best dominant. This completes the proof.

If we choose $\alpha=1$, we have
\begin{corollary}
\[
T_{n+1}^1(\beta)\subset T_n^1(2(1-\beta)\ln2+2\beta-1), \;\;\;n\in N_0.
\]
\end{corollary}

In particular for $n=0$, we have
\begin{corollary}
Let $f\in A$. Then
\[
Re f^\prime(z)>\beta\Rightarrow Re\frac{f(z)}{z}>2(1-\beta)\ln2+2\beta-1.
\]
The result is sharp.
\end{corollary}

Lastly we remark that the above result improves an earlier one due to Owa and Obradovic \cite{OO} in which they proved that if $f\in A$ satisfies $Re\;f^\prime(z)>\beta$ for $0\leq\beta<1$ and $z\in E$, then
\[
Re\frac{f(z)}{z}>\frac{1+2\beta}{3}.
\]

\vspace{10pt}

\hspace{-4mm}{\small{Received}}

\vspace{-12pt}
\ \hfill \
\begin{tabular}{c}
{\small\em  Department of Mathematics}\\
{\small\em  University of Ilorin}\\
{\small\em  Ilorin, Nigeria}\\
{\small\em {\tt ummusalamah.kob@unilorin.edu.ng}} \\
{\small\em {\tt opoola\_stc@yahoo.com}} \\
\end{tabular}


\medskip

\begin{thebibliography}{9}

\bibitem {BO}
\textsc{Babalola, K. O.} and \textsc{Opoola, T. O.}, {\it Iterated
integral transforms of Caratheodory functions and their
applications to analytic and univalent functions}, Tamkang J.
Math., {\bf37} (4) (2006), 355--366.

\bibitem {EM}
\textsc{Eenigenburg, P., Miller, S. S., Mocanu P. T.} and \textsc{Reade, M. O.}, {\it On a Briot-Bouquet differential subordination}, Rev. Roumanie Math. Pures Appl. {\bf29} (1984), 567--573.

\bibitem {TOO}
\textsc{Opoola, T. O.}, {\it On a new subclass of univalent
functions}, Matematica (Cluj) {\bf36}, 59 (2)(1994), 195--200.

\bibitem {OO}
\textsc{Owa, S.} and \textsc{Obradovic, M.}, {\it Certain subclass of
Bazilevic functions of type $\alpha$}, Int. J. Math. and Math. Sci. {\bf9} (1986), 347--359.

\bibitem {GSS}
\textsc{Salagean, G. S.}, {\it Subclasses of univalent functions},
Lecture Notes in Math. {\bf1013} (1983), 362--372.
Springer-Verlag, Berlin, Heidelberg and New York.

\end{thebibliography}
\end{document}